\documentclass[12pt]{amsproc}
\usepackage{amsmath,amssymb,amscd,stmaryrd}
\usepackage[mathscr]{euscript}

\newcommand{\sA}{\mathscr A}
\newcommand{\sB}{\mathscr B}
\newcommand{\sC}{\mathscr C}
\newcommand{\sD}{\mathscr D}
\newcommand{\sE}{\mathscr E}
\newcommand{\sJ}{\mathscr J}

\newcommand{\sT}{\mathscr T}
\newcommand{\sX}{\mathscr X}
\newcommand{\sY}{\mathscr Y}
\newcommand{\bR}{\mathbb R}
\newcommand{\bC}{\mathbb C}
\newcommand{\bg}{\mathbf g}
\newcommand{\<}{\langle}
\renewcommand{\>}{\rangle}
\newcommand{\p}{\partial}
\renewcommand{\d}{\delta}
\newcommand{\Vect}{\mbox{Vect}}

\begin{document}

\title[A Hydrodynamic Exercise in Free Probability]{A Hydrodynamic Exercise in Free Probability: Setting Up Free Euler Equations}
\author{Dan-Virgil Voiculescu}
\address{D.V. Voiculescu \\ Department of Mathematics \\ University of California at Berkeley \\ Berkeley, CA\ \ 94720-3840}
\thanks{Research supported in part by NSF Grant DMS-1665534.}
\dedicatory{Dedicated to the memory of Richard V. Kadison.}

\begin{abstract}
For the free probability analogue of Euclidean space endowed with the Gaussian measure we apply the approach of Arnold to derive Euler equations for a Lie algebra of non-commutative vector fields which preserve a certain trace. We extend the equations to vector fields satisfying non-commutative smoothness requirements. We introduce a cyclic vorticity and show that it satisfies vorticity equations and that it produces a family of conserved quantities.
\end{abstract}
\maketitle

\section{Introduction}
\label{sec1}

In $\bR^n$ equipped with the Gaussian measure, the coordinate functions can be viewed as i.i.d.\ Gaussian random variables and the divergence-free vector fields correspond to infinitesimal increments of the random variables which preserve the joint distribution. In free probability the analogue of the Gaussian random variables are the field operators on the full Fock space of $\bC^n$ with respect to vacuum expectations and the von~Neumann algebra they generate is isomorphic to the von~Neumann algebra of a free group as shown in \cite{14}. We studied in \cite{18} the infinitesimal increments for the field operators, which preserve the non-commutative joint distribution and we found that the polynomial increments among these provide a Lie algebra which is dense in the space of all such increments. Thus we have a good description of a dense part of the analogue in free probability of the divergence-free vector fields.

Since divergence-free vector fields are the basic ingredient for the hydrodynamic Euler equations we became curious about pursuing the analogy further and finding the free Euler equations. Our approach is in two steps, one formal and one analytic. The formal step consists in getting the equations on the dense Lie algebra following the general recipe in (\cite{1},\cite{2}). Of course the framework of the dense Lie algebra is much too restrictive for solutions. Fortunately, in the analytic step we are able to greatly relax the smoothness requirements on the divergence-free vector fields. Indeed, our results in \cite{18} imply that the Leray projection commutes with the number operator, which in particular makes the use of free hypercontractivity \cite{3} possible. The equations we find can be stated either in projection form (i.e., using a Leray projection) or with a pressure term. The pressure term is a cyclic gradient, roughly, and using the exact sequence for cyclic gradients we found in \cite{17}, one of the maps provides a replacement for the curl. We define in this way a cyclic vorticity, which satisfies a vorticity equation and we obtain conserved quantities which are the moments of the cyclic vorticity. In a sense this is similar to the situation on even-dimensional manifolds; however, the cyclic vorticity is not a non-commutative tensor-field, but a non-commutative scalar.

The paper has eight sections, including this introduction.

Section~\ref{sec2} contains preliminaries of a general algebraic nature from \cite{18}.

Section~\ref{sec3} deals with preliminaries concerning semicircular systems, the free probability analogue of i.i.d.\ Gaussian random variables \cite{14}, \cite{15}, \cite{16} the main aim being the results in \cite{18} about the Lie algebra of infinitesimal automorphisms of a semicircular system. We also recall some free hypercontractivity facts from \cite{3}.

Section~\ref{sec4} is an application of Arnold's procedure (\cite{1}, \cite{2}) for setting up Euler equations on a Lie algebra to the Lie algebra of infinitesimal automorphisms of a semicircular system \cite{18}. This is a formal derivation of the Euler equations, the context of non-commutative polynomials being too restrictive to expect interesting solutions.

Section~\ref{sec5} provides analytic facts which will be used to make sense of more general solutions of the free Euler equations. In particular, we deal with an algebra $\sB_{\infty,1}$ of smooth elements in the style of \cite{9}, \cite{12}, \cite{19}, which is a kind of Sobolev or Besov space. We also keep track of the special properties of the analogue of the Leray projection arising from its commutation with the number operator, or equivalently with the free Ornstein--Uhlenbeck semigroup.

In Section~\ref{sec6} we show that in view of the analytic facts about $\sB_{\infty,1}$ in the preceding section, we can make sense of the free Euler equations when non-commutative polynomials are replaced by $\sB_{\infty,1}$, roughly.

In Section~\ref{sec7} the cyclic vorticity is introduced, based on the idea that the role of the curl should be played by a map appearing in the exact sequence for cyclic gradients \cite{17}. We show that for a certain subalgebra $\sC_{\infty,1}$ instead of $\sB_{\infty,1}$ the non-commutative moments of the cyclic vorticity are conserved quantities.

Section~\ref{sec8} is devoted to concluding remarks, in particular about corresponding Navier--Stokes equations and about a possible substitute for boundary conditions.

\section{Algebraic Preliminaries (\cite{18})}
\label{sec2}

Let $s_1,\dots,s_n$ be an $n$-tuple of self-adjoint elements in a von~Neumann algebra $M$, which are algebraically free, that is there is no non-trivial algebraic relation among them. In this section we collect purely algebraic preliminaries, while later on, under additional assumptions on $s_1,\dots,s_n$ we will go beyond algebra.

The $\bC$-subalgebra of $M$ generated by $1,s_1,\dots,s_n$ will be denoted by $\bC\<s_1,\dots,s_n\>$ or $\bC_{\< n\>}$ and consists of the non-commutative polynomials in the indeterminates $s_1,\dots,s_n$.

By $\p_j: \bC_{\<n\>} \to \bC_{\<n\>} \otimes \bC_{\<n\>}$, $1 \le j \le n$, we will denote the partial free difference quotient derivations, that is the derivations so that $\p_jX_k = \d_{jk}1 \otimes 1$, $1 \le j,k \le n$. We will also use the partial cyclic derivatives $\d_j: \bC_{\<n\>} \to \bC_{\<n\>}$ defined by $\d_j = \mu\ \circ \sim \circ\ \p_j$ where $\sim$ is the flip $\sim(a \otimes b) = b \otimes a$ and $\mu: \bC_{\<n\>} \otimes \bC_{\<n\>} \to \bC_{\<n\>}$ is the multiplication map $\mu(a \otimes b) = ab$, a linear map of $\bC_{\<n\>}$ bimodules. On monomials we have
\[
\begin{aligned}
\p_j(s_{i_1}\dots s_{i_m}) &= \sum_{\{k \mid i_k=j\}} s_{i_1} \dots s_{i_{k-1}} \otimes s_{i_{k+1}} \dots s_{i_m} \\
\delta_j(s_{i_1}\dots s_{i_m}) &= \sum_{\{k \mid i_k = j\}} s_{i_{k+1}} \dots s_{i_m}s_{i_1} \dots s_{i_{k-1}}.
\end{aligned}
\]
We will also consider the cyclic gradient $\d: \bC_{\<n\>} \to (\bC_{\<n\>})^n = \bC_{\<n\>} \oplus \dots \oplus \bC_{\<n\>}$, $\d P  = \d_1P  \oplus \dots \oplus \d_nP $ and the free difference quotients gradient $\p: \bC_{\<n\>} \to (\bC_{\<n\>} \otimes \bC_{\<n\>})^n$, $\p P  = \p_1P  \oplus \dots \oplus \p_nP $.

If $\sE$ is a $\bC_{\<n\>}$-bimodule and $b \in \sE$ we consider the $\bC_{\<n\>}$-bimodules map $m_b: \bC_{\<n\>} \otimes \bC_{\<n\>} \to \sE$ so that $m_b(P \otimes Q) = PbQ$. If $b_j \in \sE$, $1 \le j \le n$ the map
\[
D_{(b_1,\dots,b_n)}: \bC_{\<n\>} \to \sE
\]
defined by $\sum_{1 \le j \le n} m_{b_j} \circ \p_j$, is a derivation of $\bC_{\<n\>}$ into $\sE$.

Since $s_1,\dots,s_n$ are algebraically free, there are evaluation homomorphisms $\varepsilon_{(a_1,\dots,a_n)}: \bC_{\<n\>} \to \sA$ where $a_1,\dots,a_n$ are elements of the unital $\bC$-algebra $\sA$, which map $s_j$ to $a_j$, $1 \le j \le n$, that is in essence $P (s_1,\dots,s_n)$ is mapped to $P (a_1,\dots,a_n)$.

If $b_1,\dots,b_n \in M$, then since $M$ is a $\bC_{\<n\>}$-bimodule, we have
\[
\frac {d}{d\varepsilon} P(s_1+\varepsilon b_1,\dots,s_n + \varepsilon b_n)|_{\varepsilon = 0} = D_{(b_1,\dots,b_n)}P
\]
where $P \in \bC_{\<n\>}$ and $\varepsilon \in \bC$. If, moreover, we also have a trace $\tau: M \to \bC$ then we have
\[
\frac {d}{d\varepsilon} \tau(P(s_1+\varepsilon b_1,\dots,s_n+\varepsilon b_n)|_{\varepsilon = 0} = \sum_{1 \le j \le n} \tau(b_j\d_jP).
\]
This means that if we endow $M^n$ with the scalar product
\[
\<(a_j)_{1 \le j \le n},(b_j)_{1 \le j \le n}\> = \sum_{1 \le j \le n} \tau(a_jb_j)
\]
then the cyclic gradient $\d P$ is the gradient at $(s_1,\dots,s_n)$ of the map
\[
M^n \ni (a_1,\dots,a_n) \to \tau(P(a_1,\dots,a_n)) \in \bC.
\]
  
If $\sE$ is a $\bC_{\<n\>}$-bimodule we shall use the notation $\Vect\ \sE$ for $\sE^u$. In particular $\sE$ can be $\bC_{\<n\>}$ or $M$.

The space $\Vect\ \bC_{\<n\>}$ is a Lie algebra under the bracket
\[
\{P,Q\} = (D_PQ_j - D_QP_j)_{1 \le j \le n}
\]
where $P = (P_j)_{1 \le j \le n}$, $Q = (Q_j)_{1 \le j \le n}$, which is the analogue of the Poisson bracket. For the analogue of hydrodynamic equations, we shall use like in Example~$5.2$ on page~20 of \cite{2}, the commutator $[P,Q] = -\{P,Q\}$.

If $\tau$ is a trace on $M$, we define
\[
\Vect\ \bC_{\<n\mid \tau\>} = \{P \in \Vect\ \bC_{\<n\>} \mid \sum_{1 \le j \le n} \tau(P_j(\d_jR)) = 0,\ \forall\ R \in \bC_{\<x\>}\}.
\]
Then $\Vect\ \bC_{\<n\mid \tau\>}$ is a Lie subalgebra of $\Vect\ \bC_{\<n\>}$. It is a non-commutative analogue of a Lie algebra of divergence-free vector fields (in an algebraic context where the vector fields can be required to be polynomial functions).

Since $s_1,\dots,s_n$ are self-adjoint, $\bC_{\<n\>}$ is a $*$-algebra. If we evaluate $P \in \bC_{\<n\>}$ at $a_1,\dots,a_n \in M$ we will have $(P(a_1,\dots,a_n))^* = P^*(a_1^*,\dots,a_n^*)$. Note also that $\d_jP^* = (\d_jP)^*$ and $\p_jP^* = \sim(\p_jP)^*$ where on $\bC_{\<n\>} \otimes \bC_{\<n\>}$ we use the involution $(\xi \otimes \eta)^* = \xi^* \otimes \eta^*$.

The involution on $\Vect\ \bC_{\<n\>}$ is defined componentwise and we have $D_{P^*}Q_j^* = (D_PQ_j)^*$. Then $P \rightsquigarrow P^*$ is a conjugate-linear automorphism of the Lie algebra $\Vect\ \bC_{\<n\>}$. In particular, the selfadjoint part
\[
\Vect\ \bC_{\<n\>}^{sa} = \{P \in \Vect\ \bC_{\<n\>} \mid P = P^*\}
\]
is a real Lie subalgebra of $\Vect\ \bC_{\<n\>}$. Similarly, since $\d R^* = (\d R)^*$ we have
\[
\Vect\  \bC_{\<n\mid \tau\>}^{sa} = \{P \in \Vect\  \bC_{\<n\mid \tau\>} \mid P = P^*\}
\]
is a real Lie algebra and
\[
\Vect\  \bC_{\<n,\tau\>}^{sa} + i \Vect\ \bC_{\<n\mid \tau\>}^{sa} = \Vect\  \bC_{\<n\mid \tau\>}^{sa}.
\]

Remark also, that if $a = (a_1,\dots,a_n) \in \Vect\  \bC_{\<n\mid \tau\>}$ and $b \in \bC_{\<n\>}$, then
\[
\tau(D_ab) = \sum_{1 \le j \le n} \tau(a_j\d_jb) = 0.
\]
In particular, if $c \in \bC_{\<n\>}$, then $\tau(D_a(bc)) = 0$ gives the ``integration by parts'' formula
\[
\tau((D_ab)c) = -\tau(b(D_ac)).
\]

\section{Semicircular Preliminaries}
\label{sec3}

From now on in this paper we shall assume $(M,\tau)$ is the von~Neumann algebra $W^*(s_1,\dots,s_n)$ generated by a semicircular system $s_1,\dots,s_n$ and $\tau$ is the unique normal trace state. So, $s_1,\dots,s_n$ have $(0,1)$ semicircle distributions and are freely independent in $(M,\tau)$, which is actually isomorphic to the $II_1$ factor generated by the regular representation of the free group $F_n$.

Let $\bC^n$ have the Hilbert space structure with orthonormal basis $e_1,\dots,e_n$ and let
\[
\sT(\bC^n) = \bigoplus_{k \ge 0} (\bC^n)^{\otimes k}
\]
be the full Fock space where $(\bC^n)^{\otimes 0} = \bC 1$ where 1 is the vacuum vector and let $l_j$, $r_j$, $1 \le j \le n$ be left and right creation operators
\[
l_j\xi = e_j \otimes \xi,\ r_j\xi = \xi \otimes e_j.
\]

Then $L^2(M,\tau)$ can be identified with $\sT(\bC^n)$ so that $s_j = l_j + l_j^*$, $1 \le j \le n$ and
\[
e_{i_1}^{\otimes k_1} \otimes \dots \otimes e_{i_p}^{\otimes k_p} = P_{k_1}(s_{i_1})\dots P_{k_p}(s_{i_p})1
\]
where $i_j \ne i_{j+1}$ $(1 \le j \le p - 1)$, $k_j > 0$ $(1 \le j \le p)$ and $P_k$, $k \ge 0$ are the Chebyshev polynomials, that is the orthogonal polynomials on $[-2,2]$ w.r.t.\ the semicircle measure (an instance of the more general Gegenbauer polynomials). Note also that the involution $\xi \to \xi^*$ on $L^2(M,\tau)$ corresponds to the antiunitary operator
\[
\sJ(ce_{k_1} \otimes \dots \otimes e_{k_m}) = {\bar c}e_{k_m} \otimes \dots \otimes e_{k_1}.
\]

The results in (\cite{18}, $7.1$--$7.5$) show that $\Vect\  \bC_{\<n\mid \tau\>}$ and $\d \bC_{\<n\>}$ have good properties as subspaces of $(L^2(M,\tau))^n \simeq (\sT(\bC^n))^n$, keeping in mind that $\d s_{i_1}\dots s_{i_p} = \d s_{i_p}s_{i_1}\dots s_{i_{p-1}}$. We have $\Vect\  \bC_{\<n\mid \tau\>} = \sum_{k \ge 0} \sX_k$, $\d \bC_{\<n\>} = \sum_{k \ge 0} \sY_k$ where $\sX_k,\sY_k \subset ((\bC^n)^{\otimes k})^n = \sX_k + \sY_k$ and
\[
\sX_k = \{((l^*_j - r^*_j)\xi)_{1 \le j \le n} \mid \xi \in (\bC^n)^{\otimes k+1}\}.
\]
Moreover $\sJ\sX_k = \sX_k$, $\sJ\sY_k = \sY_k$ and $\sX_k,\sY_k$ are orthogonal both w.r.t.\ the symmetric scalar product $\sum_{1 \le j \le n} \tau(a_jb_j)$ as well as w.r.t.\ the sesquilinear one $\sum_{1 \le j \le n} \tau(a_
jb_j^*)$.

The Hilbert spaces $\sT(\bC^n)$ and the left creation operators are part of a functor, the free analog of the Gaussian functor \cite{14}. We will need here only the free analogue of the Ornstein--Uhlenbeck semigroup $(P_t)_{t \ge 0}$ in order to use the free hypercontractivity results of \cite{3}. On $\sT(\bC^n)$ the action of $P_t$ is given by $P_t\xi = e^{-kt}\xi$ if $\xi \in (\bC^n)^{\otimes k}$. The operators $P_t$ are self-adjoint contractions and extend to all $L^p(M,\tau)$, $1 \le p \le \infty$, where they act as contractions. Moreover, on $L^{\infty}(M,\tau) \simeq M$, $P_t$ is a unit-preserving completely positive contraction.

Note also that $(P_t)^n \sX_k = \sX_k$, $(P_t)^n\sY_k = \sY_k$, where $\bC_{\<n\mid \tau\>} = \sum_{k \ge 0} \sX_k$, $\d \bC_{\<n\>} = \sum_{k \ge 0} \sY_k$, $\sX_k,\sY_k \subset ((\bC^n)^{\otimes k})^n$. In particular, if $\Pi$ is an analogue of the Leray projection, that is the orthogonal projection of $(L^2(M,\tau))^n$ onto $\overline{\bC_{\<n\mid \tau\>}}$ we have $\Pi \sJ = \sJ \Pi$, $\Pi(P_t)^n = (P_t)^n\Pi$, $\Pi\bC_{\<n\>} = \bC_{\<n\mid \tau\>}$ and $\Pi\bC^{sa}_{\<n\>} = \bC^{sa}_{\<n\mid \tau\>}$.

\section{Formal derivation of the Euler equations on the Lie algebra $\Vect\  \bC^{sa}_{\<n\mid \tau\>}$}
\label{sec4}

The recipe in I \S 4 and \S 5 of \cite{2} for Euler equations on a real Lie algebra $\bg$ endowed with a scalar product $\<\cdot,\cdot\>$ will be applied to $\bg = \Vect\  \bC^{sa}_{\<n\mid \tau\>}$ and the scalar product
\[
\<(a_j)_{1 \le j \le n},(b_j)_{1 \le j \le n}\> = \sum_{1 \le j \le n} \tau(a_jb_j).
\]
The key computation we will need to perform will provide the bilinear map $B: \bg \times \bg \to \bg$ so that
\[
\<[a,b],c\> = \<B(c,a),b\>.
\]
The scalar product being non-degenerate, $B$ is unique; however, the Lie algebra $\bg$ is not finite-dimensional and the existence of $B$ will be a consequence of the good properties of $\Pi$ the analogue of the Leray projection. The identification of $B$ is the subject of the next lemma.

\bigskip
\noindent
{\bf Lemma 1.} {\em 
Let $a,b,c \in \bC^{sa}_{\<n,\tau\>}$, $a = (a_k)_{1 \le k \le n}$, $b = (b_k)_{1 \le k \le n}$, $c = (c_k)_{1 \le k \le n}$. Then if
\[
B(c,a) = \Pi(D_ac_k + \sum_{1 \le j \le n} m_{c_j}(\sim \p_ka_j))_{1 \le k \le n}
\]
we have
\[
-\<[a,b],c\> = \<B(c,a),b\>.
\]
Moreover, we have
\[
B(a,a) = \Pi(D_a a_k)_{1 \le k \le n}.
\]
}

\bigskip
\noindent
{\em Proof.} 
We have
\[
\begin{aligned}
\<[a,b],c\> &= \sum_{1 \le j \le n} \tau((D_ab_j)c_j) - \sum_{1 \le j \le n} \tau((D_ba_j)c_j) \\
&= \sum_{1 \le j \le n} \tau(-b_j(D_ac_j)) - \sum_{1 \le j \le n} \tau((D_ba_j)c_j)
\end{aligned}
\]
where we used ``integration by parts'' since $a \in \Vect\  \bC^{sa}_{\<n\mid \tau\>}$. Note further, that
\[
\tau((D_ba_j)c_j) = \sum_{1 \le k \le n} \tau((m_{b_k}(\partial_ka_j))c_j)
\]
and note also that if $\xi,\eta \in \bC_{\<n\>}$ then
\[
\begin{aligned}
\tau((m_{b_k}(\xi \otimes \eta))c_j) &= \tau(\xi b_k\eta c_j) \\
&= \tau((\eta c_j\xi)b_k) = \tau((m_{c_j}(\sim(\xi \otimes \eta)))b_k)
\end{aligned}
\]
which gives that
\[
\tau((m_{b_k}(\partial_ka_j)c_j) = \tau((m_{c_j}(\sim \partial_ka_j))b_k)
\]
and hence
\[
\tau((D_ba_j)c_j) = \sum_{1 \le k \le n} \tau((m_{c_j}(\sim \partial_ka_j))b_k).
\]
Putting all this together gives
\[
-\<[a,b],c\> = \sum_{1 \le k \le n} \tau((-D_ac_k)b_k) - \sum_{1 \le k \le n} \tau\left(\left(\sum_{1 \le j \le n}(m_{c_j}(\sim\partial_ka_j)\right)b_k\right)
\]
where we changed indexing from $j$ to $k$ in the first sum and switched summations in the rest. This gives the formula for $B(c,a)$, the projection $\Pi$ is applied in order that $B(c,a) \in \bg$ and we get this, since
\[
\Pi\ \Vect\ \bC^{sa}_{\<n\>} = \Vect\  \bC^{sa}_{\<n\mid \tau\>}.
\]

To get the last assertion about $B(a,a)$ we must show that
\[
\sum_{1 \le j \le n} \tau((D_ba_j)a_j) = 0
\]
if $b \in \Vect\  \bC^{sa}_{\<n\mid \tau\>}$. Since $\tau$ is a trace we have
\[
\tau((D_ba_j)a_j) = \frac {1}{2} \tau(D_ba^2_j) = 0
\]
which gives the desired result.\hfill\qed

\bigskip
The hydrodynamic Euler equations for an element $v(t) \in \bg$ evolving in time is (\cite{2} top of page~20)
\[
{\dot v} = -B(v,v).
\]
In view of the Lemma we proved, the Euler equations for
\[
v(t) = (v_k(t))_{1 \le k \le n} \in \Vect\ \bC^{sa}_{\<n\mid \tau\>}
\]
will be:
\[
(\dot{v}_k)_{1\le k \le n} + \Pi(D_vv_k)_{1 \le k \le n} = 0.
\]
In our case $a - \Pi a \in \delta\bC_{\<n\>}$ so we can also figure a form of equations involving a ``pressure'' $p(t) \in \bC^{sa}_{\<n\>}$ and taking the form
\[
\dot{v}_k + D_vv_k + \delta_kp = 0\qquad 1 \le k \le n.
\]
Thus comparing with the classical equations the gradient of the pressure becomes in the free probability setting a cyclic gradient of the pressure. To continue the comparison, we did not write a continuity equation since this corresponds to the requirement $v(t) \in \Vect\ \bC^{sa}_{\<n\mid \tau\>}$, that is $v = \Pi v$.

As a final comment, the computations we did in this section did not use the assumption that $s_1,\dots,s_n$ is a semicircular $n$-tuple; however, without this assumption we have no control over $\Vect\ \bC^{sa}_{\<n\mid \tau\>}$, which may be zero and we also have no control about the projection $\Pi$.

\section{Analytic preparations}
\label{sec5}

We collect here a few analytic facts which we will use in the next section to define solutions of the free Euler equations which may not be polynomial.

Using \cite{3} an element $\xi \in L^1(M,\tau)$ can be described as a formal series $\sum_{k \ge 0} \xi_k$ where $\sum_{k \ge 0} e^{-kt} \xi_k \in L^2(M,\tau)$ (or equivalently $\sum_{k \ge 0} e^{-2kt}|\xi_k|_2^2 < \infty$) for all $t > 0$ and so that
\[
\sup_{t > 0} \left| \sum_{k \ge 0} e^{-kt}\xi_k\right|_1 < \infty.
\]
Moreover $|\xi|_1$ is precisely the above $\sup_{t > 0}$ and $\left|\xi - \sum_{k \ge 0} e^{-kt}\xi_k\right|_1 \to 0$ as $t \downarrow 0$. Note also that actually $\sum_{k \ge 0} e^{-kt}\xi_k \in M = L^{\infty}(M,\tau)$ and equals $P_t\xi$ by \cite{3}.

We may then use this to work with vector fields with $L^p$ components $1 \le p \le \infty$, that is $\Vect\ L^p(M,\tau)$, $1 \le p \le \infty$ the elements of which are $n$-tuples $\xi = \left( \sum_{k \ge 0} \xi_{k,j}\right)_{1 \le j \le n}$ where $\sum_{k \ge 0} \xi_{k,j} \in L^p(M,\tau)$, $1 \le j \le n$, where $L^p$ is viewed as a subspace of $L^1$. We shall denote by $Q_k$ the projection $Q_k\xi = (\xi_{k,j})_{1 \le j \le n}$ in $\Vect\ L^1(M,\tau)$. We have $Q_k(P_t)^n\xi = (P_t)^n Q_k\xi = e^{-kt}Q_k\xi$ and also $Q_k\xi = e^{kt}Q_k(P_t)^n\xi$.

The $L^1$-``divergence-free'' vector fields will be denoted by $\Vect(L^1(M,\tau)\mid \tau)$ and can be defined in several equivalent ways. One definition is by requiring that $\xi = (\xi_j)_{1 \le j \le n}$, $\xi_j \in L^1(M,\tau)$, $1 \le j \le n$ satisfy $\sum_{1 \le j \le n} \tau(\xi_j(\delta_jR)) = 0$ for all $R \in \bC_{\<n\>}$. Since $\d\bC_{\<n\>} = \sum_{k \ge 0} \sY_k$, $\bC_{\<n\mid \tau\>} = \sum_{k \ge 0} \sX_k$ where $\sX_k,\sY_k \subset ((\bC^n)^{\otimes k})^n$ and $(\bC^n)^{\otimes k} = \sX_k \oplus  \sY_k$ it is easily seen that $\Vect\ (L^1(M,\tau)\mid \tau)$ consists of the $\xi \in \Vect(L^1(M,\tau))$ so that $Q_k\xi \in \sX_k$. Denoting by $\Pi_k$ the projection of $\Vect\ L^1(M,\tau)$ onto $\sX_k$, this second definition of $\Vect(L^1(M,\tau)\mid \tau)$ is that $\xi \in \Vect\ L^1(M,\tau)$ satisfies $\Pi_k\xi = Q_k\xi$ for all $k \ge 0$. From here it is also easily seen that $\Vect(L^1(M,\tau) \mid \tau)$ can be defined as the closure of $\bC_{\<n\mid \tau\>}$ in $\Vect(L^1(M,\tau))$. Note also that in view of the equivalence of these 3 definitions we also have that $\xi \in \Vect(L^1(M,\tau))$ is in $\Vect(L^2(M,\tau) \mid \tau)$ iff $P_t\xi \in \Vect(L^1(M,\tau \mid \tau)$ for some $t > 0$. By $\Vect(L_{sa}^1(M,\tau))$ and $\Vect(L_{sa}^1(M,\tau)\mid \tau)$ we shall denote the corresponding subspaces of self-adjoint vector fields, that is $\xi = (\xi_j)_{1 \le j \le n}$ with $\xi_j^* = \xi_j$, $1 \le j \le n$.

We pass now to introducing a certain subalgebra of $M$ with a suitable smoothness requirement (see \cite{9}, \cite{12}, \cite{19} for related constructions). We define
\[
\sB_{\infty,1} = \{a \in M \mid [a,r_j-r^*_j] \in \sC_1\}
\]
where $a \in M$ is identified with the left multiplication operator $L_a$ on $L^2(M,\tau) \simeq \sT(\bC^n)$ and $\sC_1$ (or $\sC_1(L^2(M,\tau))$) denotes the trace-class operators on $L^2(M,\tau)$. We shall denote by $|\cdot|_p$ the $p$-norm on $L^p(M,\tau)$, $1 \le p \le \infty$ and by $\|\cdot\|_p$ the $p$-norm on the Schatten--von~Neumann classes $\sC_p$. Then
\[
\||a\|| = \|a\| + \max_{1 \le j \le n} \|[a,r_j-r^*_j]\|_1
\]
is a Banach algebra norm on $\sB_{\infty,1}$.

We will need some basic facts arising from $L^1(M,\tau)$ being the predual of $M$. There is a unique contractive linear map
\[
\Phi: \sC_1(L^2(M,\tau)) \to L^1(M,\tau)
\]
so that
\[
\mbox{Tr}(Xa) = \tau(\Phi(X)a)
\]
for all $a \in M$ and $X \in \sC_1$. If $\xi,\eta \in L^2(M,\tau)$ and $E_{\xi,\eta}$ denotes the rank one operator $\<\cdot,\sJ \eta\>\xi$, then $\Phi(E_{\xi,\eta}) = \xi\eta$. More generally every $X \in \sC_1$ can be written in the form $\sum_k E_{\xi_k,\eta_k}$ with $\sum_k |\xi_k|_2|\eta_k|_2 < \infty$ and then $\Phi(X) = \sum_k \xi_k\eta_k$. Using the $L_a,R_a$ notation for left and respectively right multiplication operators by $a$ on $L^2(M,\tau)$ we have $\Phi(R_aE_{\xi,\eta}) = \xi a\eta$ and $[R_a,X] \in \ker \Phi$ for all $a \in M$, $X \in \sC_1$. It is also easily seen that $\Phi(XL_a) = \Phi(X)a$ an $\Phi(L_aX) = a\Phi(X)$ if $a \in M$, $X \in \sC$, and $\Phi(X^*) = \Phi(X)^*$.

If  $a \in \sB_{\infty,1}$, then $[a,r_j-r^*_j] \in \sC_1$, $1 \le j \le n$ and if $b_j \in M$, $1 \le j \le n$ we have that $2^{-1} \sum_{1 \le j \le n} R_{b_j}[a,r_j-r_j^*] \in \sC_1$ and
\[
\Phi\left( 2^{-1} \sum_{1\le j \le n} R_{b_j}[a,r_j-r_j^*]\right)
\]
defines a map from $\sB_{\infty,1}$ to $L^1(M,\tau)$. This map extends the derivation $D_b$ where $b = (b_j)_{1 \le j \le n}$ from $\bC_{\<n\>}$ to $\sB_{\infty,1}$ and takes values in $L^1(M,\tau)$. Indeed, we have $[s_j,r_k-r_{k^*}] = [l_j + l_j^*,r_k-r_k^*] = 2E_{1,1} \delta_{jk}$. A straightforward computation shows that on monomials $a = s_{i_1}\dots s_{i_p}$ we have that $D_b$ and $\Phi(2^{-1}\sum_j R_{b_j}[a,r_j-r_j^*])$ are equal. We shall denote this extension of $D_b$ to a continuous linear map $\sB_{\infty,1} \to L^1(M,\tau)$ by $\sD_b$.

Consider also $\Psi: \sC_1 \to L^1(M,\tau)$ the map $\Psi(X) = \Phi(\sJ X^*\sJ)$. Since $\sJ(E_{\xi,\eta})^*\sJ = E_{\eta,\xi}$ we will have $\Psi(E_{\xi,\eta}) = \eta\xi$. Let ${\tilde \delta}: \sB_{\infty,1} \to L^1(M,\tau)$ be defined by ${\tilde \delta}_j(a) = \Psi(2^{-1}[a,r_j-r_j^*])$, $1 \le j \le n$. If $a \in \bC_{\<n\>}$ we have ${\tilde \delta}_j(a) = \delta_j(a)$, so that ${\tilde \delta}_j$ is an extension of $\delta_j$ to $\sB_{\infty,1}$.

An easy computation also proves the following.

\bigskip
\noindent
{\bf Lemma.} {\em 
$\sD_b: \sB_{\infty,1} \to L^1(M,\tau)$ is a derivation and then $\sD_{b^*}(a^*) = (\sD_b(a))^*$, in particular if $b=b^*$ then $\sD_b(a^*)=(\sD_b(a))^*$. Moreover we have
\[
|\sD_ba|_1 \le 2^{-1} \sum_{1 \le j \le n} \|b_j\|\|[a,r_j-r_j^*]\|_1.
\]
Similarly we have
\[
|{\tilde \delta}_j(a)|_1 \le 2^{-1}\|[a,r_j-r_j^*]\|_1.
\]
}

\section{The free Euler equations in $\sB_{\infty,1}$}
\label{sec6}

In this section we show that the free Euler equations we found in the formal setting of section~\ref{sec4}, with some adjustments, make sense also for $\Vect(\sB_{\infty,1}^{sa} \mid \tau)$, that is for divergence-free vector fields with components in $\sB_{\infty,1}^{sa}$.

We recall the equations were
\[
({\dot v}_k)_{1 \le k \le n} + \Pi(D_vv_k)_{1 \le k \le n} = 0
\]
and the requirement $v(t) \in \Vect\ \bC^{sa}_{\<n|\tau\>}$ was playing the role of the continuity equation. We shall assume that each $v_k(t) \in \sB^{sa}_{\infty,1}$, $1 \le k \le n$ where $t \in [0,T)$ some $T > 0$ is differentiable as a function of $t$ with values in $L^1(M,\tau)$ (which is a weaker requirement than as a function with values in the Banach space $\sB_{\infty,1}$). We shall also assume that 
\[
v(t) = (v_k(t))_{1 \le k \le n} \in \Vect(L^1(M,\tau)\mid \tau)
\]
for all $t \in [0,T]$.

With the above requirements we get clearly that $({\dot v}_k(t))_{1 \le k \le n} \in \Vect(L_{sa}^1(M,\tau)\mid \tau)$.

Since $(v_k(t))_{1 \le k \le n} \in \sB_{\infty,1}$ we have that $\sD_{v(t)} v_k(t) \in L^1(M,\tau)$ by the results of section~\ref{sec5}. This suggests that we replace $D_v(t)$ which we had defined only to act on $\bC_{\<n\>}$, by $\sD_{v(t)}$ which acts on $B_{\infty,1}$.

The resulting $\sD_{v(t)}v(t)$ being in $\Vect\ L^1(M,\tau)$ we face the problem that the free Leray projection is defined in $\Vect\ L^2(M,\tau)$. Since the Leray projection commutes with the number operator we can use instead of $\Pi$ the projections $\Pi_m$, $m \ge 0$ and we get a sequence of equations
\[
-Q_m{\dot v}(t) = \Pi_m(\sD_{v(t)}v_k(t))_{1 \le k \le n},\ m > 0.
\]
An essentially equivalent way is to use the Ornstein--Uhlenbeck semigroup with some $\epsilon > 0$. The equations become
\[
-P_{\epsilon}{\dot v}(t) = \Pi P_{\epsilon}(\sD_{v(t)}v_k(t))_{1 \le k \le n}.
\]

Other possibilities one may consider would be to use a suitably defined cyclic gradient of pressure term. The equations give
\[
-{\dot v}(t) = (\sD_{v(k)}v_k(t))_{1 \le k \le n} + q(t)
\]
where
\[
q(t) \in \Vect\ L^1_{sa}(M,\tau)
\]
is continuous as a function of $t$ and should satisfy $\Pi_kq(t) = 0$ for all $k \ge 0$. This implies each $Q_kq(t) \in \delta\bC_{\<n\>}$. In view of the exact sequence for cyclic gradients \cite{17}, this gives
\[
\sum_{1 \le j \le n} [s_j,Q_kq(t)] = 0.
\]
It is easy then to infer from here that
\[
\sum_{1 \le j \le n} [s_j,q(t)] = 0.
\]

\section{Cyclic vorticity and conserved quantities}
\label{sec7}

In the classical setting applying the curl to the Euler equations makes the gradient of the pressure disappear and leads to the vorticity equations and the expectation values of powers of the vorticity, which is a differential form, are roughly the source of conserved quantities. In the free probability setting, the pressure term is a cyclic gradient and instead of the curl we have a map in the exact sequence we gave in \cite{17} which will make the pressure term disappear. The new quantity we get is in $M$ and expectation values of its powers, or differentiable functions of it, provide conserved quantities, a situation somewhat reminiscent of the classical Euler equations in even dimension. Since the analogy works reasonably well, we will call the quantity we get, the cyclic vorticity.

If $v \in \Vect(L_{sa}^1(M,\tau))$ we define
\[
\Omega = i \sum_{1 \le j \le m} [s_j,v_j]
\]
to be its cyclic vorticity. This is just $i\theta(v)$ with $\theta$ the extension to $L^1$ of the map in the exact sequence for cyclic gradients \cite{17} and the coefficient $i$ has been added to make sure that
\[
\Omega \in \Vect(L^1_{sa}(M,\tau)).
\]

Assume now $v(t) \in \Vect(\sB_{\infty,1}^{sa}) \cap \Vect(L_{sa}^1(M,\tau)\mid \tau)$ is differentiable as a $\Vect(L^1(M,\tau))$-valued function of $t \in [0,T)$ and satisfies the free Euler equations in the form with pressure-term, discussed at the end of section~\ref{sec6}, that is
\[
-{\dot v}(t) = (\sD_{v(t)}v_k(t))_{1 \le k \le n} + q(t)
\]
where $q(t) \in \Vect(L_{sa}^1(M,\tau))$ is continuous as a function of $t$.

Since $\sum_{1 \le j \le n}[s_j,q_j(t)] = 0$ as discussed in section~\ref{sec6}, we have
\[
\begin{aligned}
i{\dot \Omega}(t) &= -\sum_{1 \le j \le n} [s_j,{\dot v}_j(t)] \\
&= \sum_{1 \le j \le n}([s_j,\sD_{v(t)}v_j(t)] + [s_j,q_j(t)]) \\
&= \sum_{1 \le j \le n}(\sD_{v(t)}[s_j,v_j(t)] - [\sD_{v(t)}s_j,v_j(t)]) \\
&= \sum_{1 \le j \le n}(\sD_{v(t)}[s_j,v_j(t)] - [v_j(t),v_j(t)]) \\
&= i^{-1}\sD_{v(t)}\Omega(t). 
\end{aligned}
\]
This gives the cyclic vorticity differential equation
\[
-{\dot \Omega}(t) = \sD_{v(t)}\Omega(t).
\]
Note also that $\Omega(t) \in \Vect(\sB_{\infty,1}^{sa})$.

The usual vorticity is divergence-free and this provides a further equation. To see what the analogue of this is in our setting, we need to return to the exact sequence for cyclic gradients \cite{17} and the range of the map $\theta$. The range of $\theta$ is contained in $\mbox{Ker } C$, where $C$ is the map $Cs_{i_1}\dots s_{i_p} = \sum_{1 \le j \le p} s_{i_j}\dots s_{i_p}s_{i_1}\dots s_{i_{j-1}}$ also in $\mbox{Ker } \tau$. On the other hand $Ca = \sum_{1 \le j \le n} s_j\delta_ja$. Using the discussion at the end of section~\ref{sec5} and the Lemma, the map $C$ has an extension ${\tilde C}a = \sum_{1 \le j \le n} s_j{\tilde \delta}_ja$. Introducing the subalgebra $\sC_{\infty,1}$, the closure of $\bC_{\<n\>}$ in $\sB_{\infty,1}$, we find that if $v(t) \in \sC_{\infty,1}$, then we have
\[
{\tilde C}\Omega(t) = 0 \mbox{ and } \tau(\Omega(t)) = 0.
\]

As a consequence of this differential equation, if $m \in {\mathbb N}$ we have
\[
-\frac {\partial}{\partial t} \tau(\Omega^m(t)) = \tau(\sD_{v(t)}\Omega^m(t)).
\]
If $\Omega^m(t)$ would be in $\bC_{\<m\>}$, the fact that $v(t) \in \Vect((L^1(M,\tau)\mid \tau)$ would imply $\tau(\sD_{v(t)}\Omega^m(t)) = 0$ and hence $\tau(\Omega^m(t))$ would be a conserved quantity. Building on this observation, we will use $\sC_{\infty,1}$ the subalgebra of $\sB_{\infty,1}$ which is the closure of $\bC_{\<n\>}$ with respect to the norm $\||\cdot\||$. Using the lemma at the end of section~\ref{sec5}, we have

\bigskip
\noindent
{\bf Lemma.} {\em 
If $v \in$ {\rm Vect}$(M,\tau)$ and $a \in \sC_{\infty,1}$, then we have
\[
\tau(\sD_va) = 0.
\]
}

We then conclude that

\bigskip
\noindent
{\bf Theorem.} {\em 
If $v(t) \in$ {\rm Vect}$(\sC^{sa}_{\infty,1}) \cap$ {\rm Vect}$(L_{sa}^1(M,\tau) \mid \tau)$ is differentiable as a {\rm Vect}$(L^1(M,\tau))$-valued function of $t \in [0,T)$ and satisfies the free Euler equations in the form with pressure-term, then $\tau(\Omega^m(t))$ is a constant where $\Omega(t)$ is the cyclic vorticity and $m > 0$.
}

\bigskip
Clearly the theorem can be stated more generally, that $\tau(\varphi(\Omega(t)))$ is constant, where $\varphi$ is a bounded Borel function, since we know that $\Omega(t)$ is bounded and its distribution is constant.

\section{Concluding remarks}
\label{sec8}

We collect here various remarks related to the free Euler equations.

What is the analogue of boundary conditions for our free Euler equations? A possible answer may be: the requirement that the non-commutative vector field $v(t)$ be the generator of a one-parameter group of automorphisms.

We didn't go beyond Euler equations in this paper, but clearly free Navier--Stokes equations can be imagined where the viscosity term would involve the number operator. In particular, the commutation of the Leray projection with the number operator can still be expected to help.

We should also mention that replacing Gaussian random variables by non-commutative generalizations, was considered in other questions concerning fluids, quite early, in the pioneering paper \cite{5} which led to later non-commutative probability work \cite{6}, \cite{7}.

We should mention that transportation for classical PDE \cite{8}, \cite{13} being relevant to Euler equations and having now free transportation developments \cite{4}, \cite{10}, \cite{11}, in particular the last two paper achieving very strong results, there seem to be reasons for optimism about going deeper into the analytic problems of the free Euler equations.

\end{document}